# FROM PRACTICAL GEOMETRY TO THE LABORATORY METHOD: THE SEARCH FOR AN ALTERNATIVE TO EUCLID IN THE HISTORY OF TEACHING GEOMETRY

Marta Menghini

Dipartimento di Matematica, Sapienza Università di Roma

marta.menghini@uniroma1.it

*This paper wants to show how practical geometry, created to give a concrete help to people involved in trade, in land-surveying and even in astronomy, underwent a transformation that underlined its didactical value and turned it first into a way of teaching via problem solving and then into an experimental-intuitive teaching that could be an alternative to the deductive-rational teaching of geometry. This evolution will be highlighted using textbooks that proposed alternative presentations of geometry.*

## INTRODUCTION

As is well known, *practical geometry* was present from the very beginning of the history of geometry. Different populations, such as the Egyptians, the Babylonians, the Chinese and the Indians developed practical geometric skills. A practical geometry was present in the Greek world too; however, "basing their work on the practical geometry of the Egyptians", the Ancient Greeks developed "a system of logic which culminated with the great work of Euclid" (Walker Stamper, 1909, p. 4). As we know, the *logical-deductive* aspect of Euclid's *Elements* greatly influenced the teaching of geometry.

The Romans were mainly interested in practical geometry (land-surveying and the engineering of warfare), and in its *teaching in their schools*. In this environment, we find one of the first *didactical appreciations* of practical geometry: "Leaving its use in warfare, many maintain that this science, different from others, is not useful when it is part of the knowledge, but *at the moment in which it is learned* (Quintilianus, *Istitutio Oratoria*, I, 34 etc.). With these words *Quintilianus* attributed to practical geometry, through the attention to the *process* of learning, the same educational role that *proving* has within rational geometry.

Many attempts have been made to teach geometry following a path different from that proposed in the Elements of Euclid (Barbin & Menghini, 2013) The most diffused way was to change the order of the theorems and of the problems, while accepting hypothetical constructions and a limited use of arithmetic. In this paper we will look at those attempts whose main purpose was not in a (new) logical order for geometric deduction, but rather in an "exploration" of geometry that kept a strong link with its practical origin, and whose contents and methods can be seen as an evolution of the medieval practical geometry.

# 1. "*DE PRACTICA GEOMETRIAE*"

We will begin our excursus among the *textbooks* of practical geometry with Fibonacci's *De practica geometriae* (1223). This work by Leonardo of Pisa (Fibonacci) gave rise to that stream of practical geometry which characterised the geometry of the Late Middle Ages and presented not only the very practical rules and methods for calculating distances and areas for land-surveyors, but also problems that, towards the end of the Middle Ages, turned into "*mathematical games*" (as in Leon Battista Alberti's *Ludi matematici*).

Fibonacci's text was largely used in the last grades of the Italian *Scuole d'abaco*, which were parish-schools for pupils who wanted to learn a trade. In these schools pupils had to memorize only the rules, but in the text of Fibonacci there was something that went beyond this idea.

Fibonacci didn't write his book in contrast the Euclidean one. On the contrary, he often refers to Euclid's propositions. His geometry is simply a different thing. Fibonacci never speaks of axioms or theorems; at the beginning he lists definitions and "principles", which broadly correspond to axioms or to constructions that can be done. He uses *numbers*, arithmetic, and practical examples; proofs are often only verifications with numbers. In the case of the theorems from the 2nd book of Euclid, concerning geometric algebra, he reports the proofs "translated" in an algebraic language.

Fibonacci doesn't present drawing or measuring instruments. He doesn't even refer to angles. His geometry has the scope of "measuring all kinds of fields" and "dividing fields among partners". In fact one of the first rules shown by Fibonacci concerns the *calculation of the area of a square*. As he always does, the rule is given with an example (p. 14)[1].

> Given a quadrilateral, equilateral and equiangular field having 2 rods on each side, I say that its area is to be found by multiplying side *ac* by its adjacent side *ab*, namely two rods by two rods.
>
> 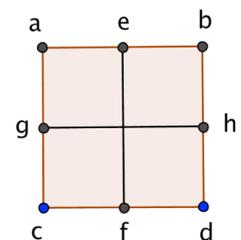
>
> Let lines *ab* and *cd* be divided into two equal parts [this is one of the allowed constructions listed at the beginning] at points *e* and *f*, and draw the line *ef*. Likewise […] draw line *gf*. Thus quadrilateral *abcd* has been divided into four perpendicular squares, each of which is measured by one rod on a side. Thus there are 4 plane rods in the entire square quadrilateral *abcd*.

This explanation will be sufficient for all later practical geometry. But Fibonacci wants to show that in fact the four quadrilaterals are squares. The "proof" consists in an *observation*: being *ef* equal to and equidistant from *ac* and *bd*, and being *gh* equal to and equidistant from *ab* and *cd,* also *ae* is equidistant and equal to *cf.* Therefore the angle *aef* is equal to the right angle *acf*, a.s.o.

---

[1] For the English translation, as well as for the pages for the quotes, we refer to the edition of Hughes (2008).

Fibonacci uses the theorem of Pythagoras to measure areas, and he proves it by referring to Euclid's theorems concerning the similarity of the involved triangles; in the triangle *agd*, we have (p. 69)

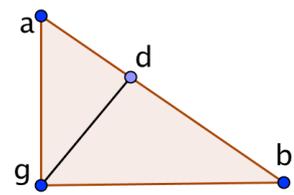

*db:bg = gb:ab*, hence *db·ba = bg²*. Similarly *ad·ab = ag²*.

So *db·ba + ad·ab = bg² + ag²*. But the first part equals the square on *ab*, as proved by Euclid.

Then Fibonacci gives the rule for the area of a triangle: "To find the square measure, namely the areas of all triangles, multiply half the cathete [the height] by the whole base or half the base by the whole catheter" (p.70). Fibonacci shows, only in the case of a rectangular triangle, that a triangle can be seen as half a rectangle. For the other triangles, he tries to come back to right-angled triangles by dividing the figure by using the height. He gives many examples in order to consider various positions of the height.

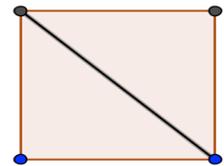

The many problems concerning "dividing fields" into equivalent parts are more abstract, although they originate from concrete problems. Here we start to find "puzzles", challenges. For instance "To divide a triangle in two equal parts by a line drawn from a given point on a side" (p. 189):

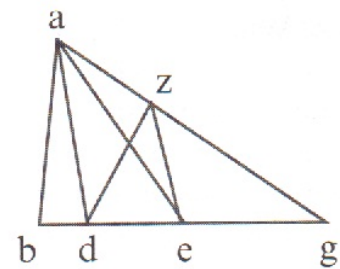

In the triangle *abg*, consider the point *d*. I will divide side *bg* in two equal parts at point *e*, and I will join lines *ad* and *ae*. And through point *e* I will draw line *ez* equidistant from line *ad*, and I will draw line *dz*. I say therefore that triangle *abg* is divided in two by line *dz*. The proof follows. Two equal triangles *ade* and *adz* sitting on base *ad* with sides *ad* and *ez* are equidistant from one another. To both triangles add triangle *abd*. The two triangles *abd* and *adz* become quadrilateral *abdz* equal to two triangles *abd* and *ade* that is triangle *abe*. But triangle *abe* is half of triangle *abg*. Whence quadrilateral *abdz* is half of triangle *abg*. What is left, namely triangle *adg*, is the other half of triangle *abg*. Therefore triangle *abg* is divided into equal parts at point *d* by line *dz*.

Numerically: *ed·ga/gd = az*, for instance if *ab*=13, *bg* = 14, *ga* = 15, then *gz* = 11 e 2/3.[2]

These kinds of problems are not found in Euclid. In the following centuries, they became typical problems in texts of practical geometry.

In the chapter about "measuring heights, depths, and longitudes of planets", we find practical problems as (p. 346)

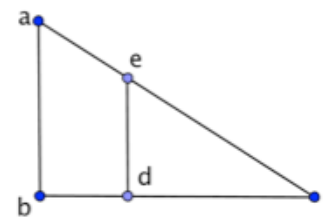

If you wish to measure a height [*ab*], fix a staff [*ed*] perpendicular to the ground. Step away from the staff and the object you wish to measure. Stoop down the ground level from where you can see the top of the object across the top of the staff, and mark the place from where you looked [*c*].

---

[2] The triangle 13 – 14 – 15 was often used in practical problems, as it has a "good" area (42).

Considering the involved similar triangles, Fibonacci explains how to find the height. He gives different examples, always repeating how the rules are used ("This can be shown with numbers. Let *ed* be 5 palms [...]").

With Fibonacci's text the students' learning is facilitated because the rules are introduced through numerical examples, and theorems are "proved" by substituting numbers to the considered lengths, and, sometimes, with the help of algebra (as the theorem of Pythagoras). On the other hand, the student is challenged with problems, such as the division of fields.

Books similar to Fibonacci's, may be with a little more Euclidean proofs, lasted for over 300 years and influenced the teaching not only in the old universities, but also in the first secondary schools of the 16[th] century. Examples are the texts of Luca Pacioli in Italy (1494) and of Orontius Finaeus in France (1556).

## 2. "GEOMETRIAE LIBRI": THE WAY TO GEOMETRY[3]

In Renaissance, Italy lost its mathematical "supremacy" to France. In fact, in the 16[th] century we find the interesting work of Pierre de la Ramée (Petrus Ramus). The text by Ramus (we refer to the edition of 1569) should be about practical geometry. But, even though Ramus presents more drawing and work instruments then Fibonacci, he is also a philosopher and an educator and explicitly criticizes the presentation of Euclid.

Ramus writes in the preface that "geometry is the art of measuring well" (*geometria est ars bene metiendi*). To measure well it is necessary to consider the nature and "affections" of everything that is to be measured: to compare such things one with another, to understand their reason, proportion and similarity ... This aim of Geometry appears much more beautiful when one observes "astronomers, geographers, land-meaters, sea-men, enginers, architects, carpenters, painters, and carvers" in their description and measurement of the "starres, countries, lands, engins, seas, buildings pictures, and statues or images".

Ramus starts with an extended *presentation* of geometric entities and figures through drawings, measures and simple constructions. We could call this an *observational* geometry, through which the reader becomes acquainted with the figures and their properties. We can say that in Ramus we find the concept of the different *levels* of geometric understanding. To reach his aim, Ramus also presents *unusual* objects and relations. Some examples (Book 4):

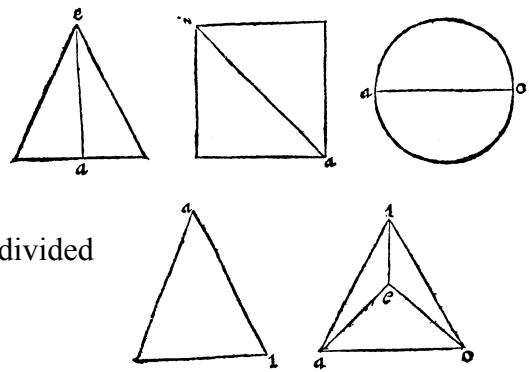

6. [The Diameter is a right line inscribed within the figure by its center]. The diameters in the same figure are infinite.

- if a figure has all equal diameters it is a circle.

11. A prime or first figure, is a figure which cannot be divided into any other figures more simple then it selfe.

---

[3] *The way to geometry* is the title of the English translation of Bedwell (Ramus, 1636) that we use here for the English quotes.

Ramus introduces drawing instruments (ruler and compass), then, also referring to previously mentioned properties he states

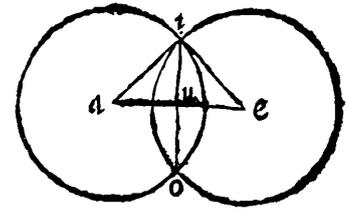

12. If two equall peripheries from the ends of a right line given, doe meete on each side of the same, a right line drawne from those meetings, shall divide the right line given into two equall parts (Book 5).

Ramus explains that the segment *ae* of the figure is divided into two equal parts by *io*, using the equality of triangles. Then he shows how to construct a triangle given three sides.

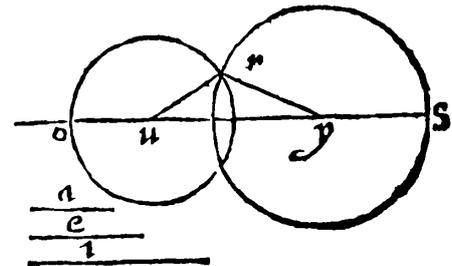

10. If of three right lines given, any two of them be greater than the other, and peripheries described upon the ends of the one, at the distances of the other two, shall meete, the rayes from that meeting unto the said ends, shall make a triangle of the lines given (Book 6).

In all cases Ramus explains his constructions, he proves that they work, linking the theoretical to the practical aspect. The text is full of drawings and of properties of the figures that we normally don't find in textbooks. He brings the reader to observe in an active way; and he supplies the proofs. We also find problems about the partition of fields as in Fibonacci:

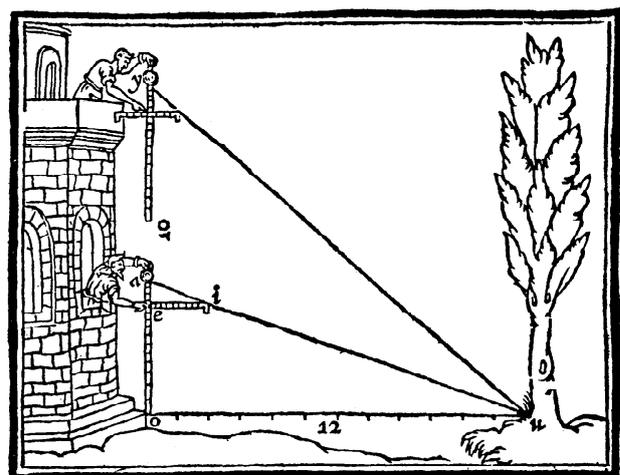

10. If a right line be drawne from the toppe of a triangle, unto a point given in the base (so it be not in the middest of it) and a parallell be drawne from the middest of the base unto the side, a right line drawne from the toppe of the sayd parallell unto the sayd point, shall cut the triangle into two equall parts (Book 7)

We have omitted the proof, which is the same as in Fibonacci.

Then Ramus arrives to the "classical" problems of practical geometry, giving rules (and instruments) for the measuring of segments (Book 9, the quotation refers to the lower part of the figure):

7. If the sight be from the beginning of the Index right or plumbe unto the length, and unto the farther end of the same, as the segment of the Index is, unto the segment of the transome, so is the heighth of the measurer unto the length.

Let therefore the segment of the Index, from the toppe, I meane, unto the transome be 6. parts. The segment of the transome, to wit, from the Index unto the opticke line be 18. The Index, which here is the heighth of the measurer, 4. foote: The length, by the rule of three, shall be 12. foote. The figure is thus, for as *ae*, is to *ei*, so is *ao*, unto *ou*, […] for they are like triangles.

Up to this point, Ramus' text is concerned with the *measuring of lengths*. Only after a long chapter on quadrilaterals and their properties, the way to measure the *area of a rectangle* by dividing it in unit squares is shown.

In Ramus' text measurement is still the main aim of geometry, but the activity of measuring is preceded by a very interesting and original part of observation, a "play" with figures. The theorems that are proved mainly belong to the classical tradition, but their proofs are strongly supported by geometric construction, and the use of drawing is added to the practical activity of measuring by means of proper instruments.

## 3. "ELEMENTS DE GÉOMETRIE": GEOMETRIC CONSTRUCTIONS

In 1741, again in France, Alexis Clairaut wrote his Elements de Géométrie. His first chapter is about the *measurement of fields*, nevertheless Clairaut was *not* interested in teaching a practical geometry. With Clairaut we see a shift from measurement as a goal to measurement as *a means* to teach geometry via problems. This is seen by the fact that the part about measurement *doesn't contain numbers*; there is only a hint at the necessity for a comparison with a known measure.

The aim of Clairaut is to solve a problem "constructing" the elements he wants to measure. The focus is on the *process* of constructing and in a *narrative* method.

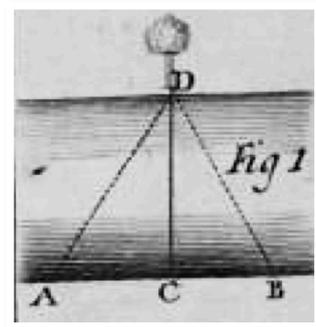

A person placed on *D*, on the bank of a river, wishes to ascertain how far it is from the other bank *AB*. It is clear, in this case, that to have this measure, it is necessary to take the shortest of all the right lines *DA*, *DB*, etc., which may be drawn from the point *D* to the right line *AB*. Now it is easy to observe that the line wanted is DC, which must not lean to the line *AB* more in the direction of the right hand than of the left. On this therefore, which is called a perpendicular, we must move our known measure [...] We are therefore under the necessity of finding a method of drawing perpendiculars (p. 2).

Clairaut shows other cases in which a perpendicular is needed, for instance, to draw a rectangle. And then he goes on with the construction:

The point C might be found by repeated trials, but this method would leave the mind unsatisfied [...] Take a common measure, a rope line, or a pair of compass with a certain opening, according as you may going to operate either on the ground or on paper....

Even if a more detailed part has been omitted, this quotation shows the precision of the description. There is no definition of perpendicular, but the term is clear when it appears. We can note here the importance of the "moment" of learning, mentioned by Quintilianus.

Clairaut's text doesn't contain proofs, but constructions and arguments. The third part of the text also presents concrete models for space geometry.

Numbers appear rarely in Clairaut's treatise: mainly to show that the area of a rectangle is given by the product of its base by the height (dividing the rectangle in horizontal stripes whose height is one, and then each stripe into unit squares). The area of a triangle is then half the product of its base and height, because "we may easily perceive that this figure [the rectangle] transversely divided by the line *AC*, which is called a diagonal, resolves itself into two equal triangles" (p. 12). Clairaut seems not to be worried, as Fibonacci was, of showing that the property holds for all triangles.

The second part of the text is more *rigorous*, as Clairaut himself asserts, not because there are proofs, but because the allowed instruments are only a ruler and a compass. Clairaut continues via his problem solving method to transform, for example, a rectangle into another rectangle with equal area.

> Let it now be proposed to change the rectangle *ABCD* into another, which shall have the same surface but with the altitude *BF*. [...] the new rectangle (the altitude of which is to be greater than *BC*) must necessarily have a base less than *AB*; [...] All therefore we have to do is to divide the line *AB* in such manner that *AB* shall be to *GB*, as *BF* to *BC*. This (according to Part I, XLI) will be done by drawing the line *FA*, and, from the point *C*, the line *GC* parallel to *FA* (p. 68).

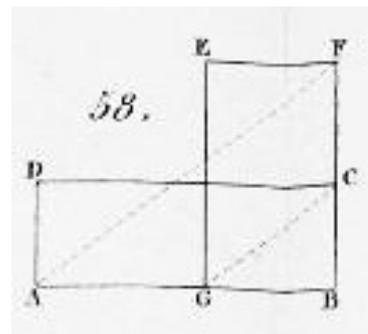

Clairaut introduces proportions and gives numerical examples.

To sum two squares "*we can transform the smallest square into a rectangle whose height is the side of the other square, so obtaining a new rectangle*".

> To sum two squares in order to obtain another square, we can work on the decomposition of the two squares, as in following figure, where *H* is taken so that *DH* = *CF* (and consequently *HF* = *fh* = *DC*) (p.76).

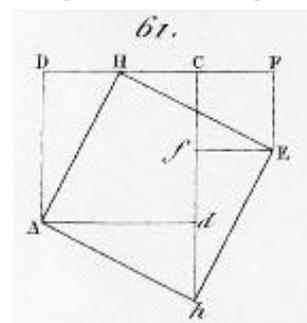

So the theorem of Pythagoras comes to mind when speaking of equivalent figures and of comparison of figures (though Clairaut doesn't explicitly mention Pythagoras).

Clairaut doesn't make use of the classical theorems of Euclid, but he does prove the theorems concerning the angles in a circumference. He then goes back to problems of land-surveying; this has the purpose of using properly the protractor so as to construct a similar figure on the paper "Knowing the distances between the given points *ABC*, we want to find their distances to a fourth point *D* from which the three are seen" (Part 3, n. 91). Then two circumferences can be constructed, which have as a periphery angle the angle *adb* and the angle *bdc*....

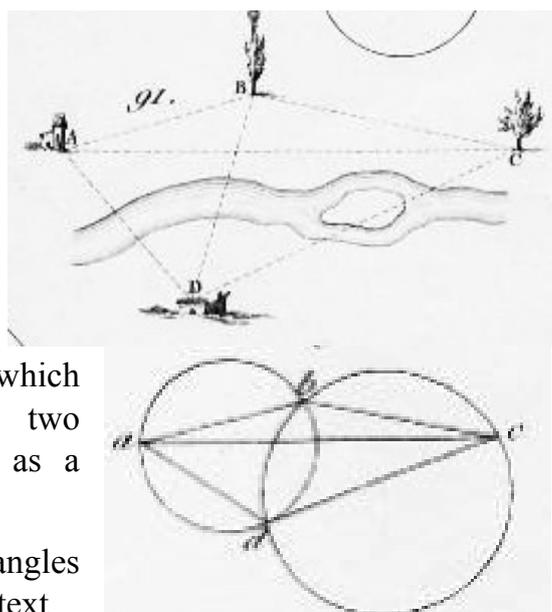

This application of the theorem about the periphery angles is one of the few challenging problems of Clairaut's text.

The challenges of Clairaut are not in solving difficult problems, but rather in the request for the learner to follow him in the process of the solution to a problem, by paying attention to the needed construction.

Clairaut's success came a century after his text was written. In 1836 it was translated for use in the Irish national schools (Clairaut, 1836)[4] and was reprinted in French in 1852 and officially adopted. It was also used in Italian Technical Schools (which correspond to the first three years of technical instruction) till to the beginning of the 20[th] century.

## 4. "ANFANGSGRÜNDE DER GEOMETRIE": DRAWING, CONSTRUCTIONS, PROOFS

In 1846 Franz (Ritter von) Mocnik wrote a text that had many reprints and adaptations up to the beginning of the 20[th] century. It was used in the Austro-Hungarian territories, including the northern of Italy before the Italian unification, and was mainly addressed to technical and professional schools (see Mocnik, 1873). It is another text for the practical application of geometry, but it is now a textbook aimed at an extended school-system; and as a result pedagogical theories start to emerge.

Mocnik divides his text into a first part called *Formenlehre*, which includes free-hand drawing, some hints about projective geometry and practical measurements, and a second part called *Grundlehre*, which includes constructions with ruler and compass (as in Clairaut) and also proofs. *Formenlehre* reminds us of the educator Johann Heinrich Pestalozzi, who was particularly successful in Mittel- and North-European countries, and who supported an initial approach to teaching based on intuition. Projective geometry is included in *Formenlehre* because it describes the way we see.

Mocnik starts with practical indications on how to represent points, on how to draw a line by free-hand or with a ruler. In fact, he gives many exercises that require free-hand drawing. There are also problems applied to land-surveying, as

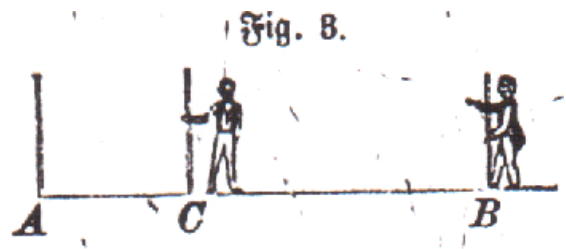

If we want to plant a pole, in a field, that must be aligned with two other poles *A* and *B*, we need to be behind pole *B* while an assistant goes with pole *C* to where the pole should be planted. Then we can indicate to him with our hand where he should move until we can see that pole *C* is aligned with *A* and *B* (p. 9, translated by the author).

The concept of distance is primitive, and an there is an exercise that asks, for instance, to draw five horizontal lines at equal distances; then Mocnik suggests to compare approximately the length of two segments and to draw their sum or difference (there is no reference to numbers). Also approximate multiples or submultiples of a segment have to be found through the drawing of them.

---

[4] This is the translation that we used for the quotes of the first two parts.

Some information is given about new and old unit measures, and about measuring instruments. The angle is introduced as the difference between the inclination of two straight lines, and the perpendicular [to a line *AB*] is defined as in Clairaut, as a line *which must not lean to the line AB more in the direction of the right hand than of the left.* But Mocnik also adds that it will form two equal angles with the line *AB*.

This geometric description and others (such as the diagonals of a square) help pupils to understand the concept and thus be able to draw perpendiculars free-hand. In a field a perpendicular has to be drawn using a special instrument.

Measuring instruments are only introduced after that the concept has been explored through free-hand drawing. But the use of instruments is not an aim, their use is not often suggested. For instance, after the introduction of the protractor, Mocnik shows how to draw a 60° angle (p. 30)

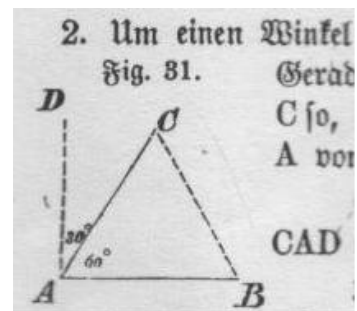

> We take a segment *AB*, and a point *C* whose distance from *A* and *B* should be as the distance from *A* to *B*. Then we draw *AC*. *BAC* = 60°. If we then draw the perpendicular *AD* to *AB*, then *CAD* = 30° (translated by the author).

Mocnik also gives a method to "construct" our own protractor for angles less then 16 degrees (approximating the values of the tangent), but he also suggest to guess the value of an angle.

Mocnik doesn't explain all the rules and methods that he gives, nevertheless we cannot say that his practical geometry is in the "medieval" meaning of the term. The suggested work aims at giving the pupil a familiarity with the topics. In a certain sense it has something in common with Ramus'. And in fact Mocnik slowly increases the level of rigor.

He describes solids in the space and gives the first elements of perspective, for instance:

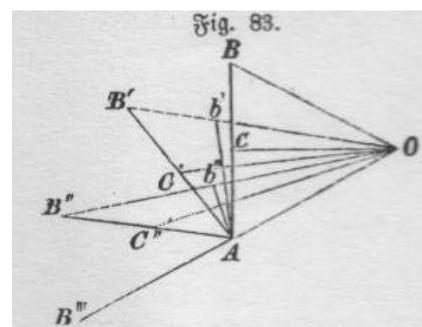

> If a line [a segment] is rotated from the position *AB* ... around *A* to the positions *AB'*, *AB''*, *AB'''*, it can be seen under different angles and with different lengths. The apparent length of the segments is obtained when, from point *A*, we draw – inside the corresponding angle of view – a perpendicular to the central ray (p. 64, translated by the author).[5]

Congruent triangles are, in this first part of the text, simply triangles that have three equal sides and three equal angles. Mocnik suggests drawing them using *motions*. In the part concerning *Grundlehre*, the first proof is about the equality of two triangles that have three equal sides. Indications are given on how to construct, with ruler and compass, a triangle that is equal to another one, and then how to construct a triangle given three elements.

---

[5] The drawing doesn't seem very precise, but the description is correct.

Theorems are always followed by exercises that require constructions. We find practical applications of the congruence of triangles (p. 93):

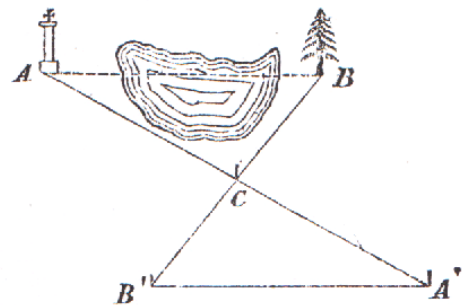

Fig. 119.

> To determine the distance of two points in a field, if this cannot be measured directly due to an obstacle between them, but it is possible to measure the distance from a third point to both…

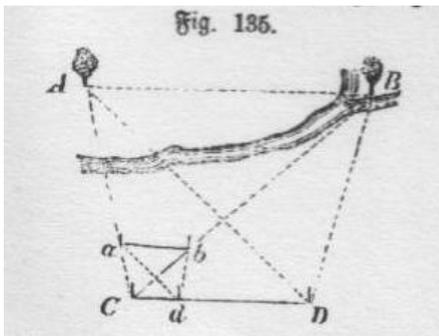

Fig. 135.

...and practical applications of the similarity (p.104): *"To determine the distance between two un-reachable points"* you have to fix two points *C* and *D* and measure the angles *CDB*...

Only at this point Mocnik introduces the measurement of areas. The area of a square is determined as in Clairaut; Mocnik at first sees the figure as subdivided in rectangles and then in unit squares. He proposes many numerical examples of this rule.

The theorem of Pythagoras too is proved in the way of Clairaut, but Mocnik also proposes to draw the figure used for the proof (see section 3) on paper and then to cut its various parts in order to re-obtain the considered squares.

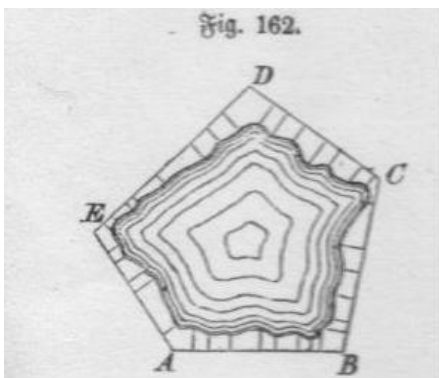

Fig. 162.

Again we find practical problems, such as to represent parts of the land on paper. Here Mocnik also describes instruments used by surveyors.

And again we come back to theorems. Mocnik proves those concerning the circumference and also conic sections and geometry in space.

## 5. "LEHRBUCH DER ELEMENTAR GEOMETRIE": GEOMETRIC TRANSFORMATIONS

In the second half of the 19[th] century many texts tried to modify the introduction to geometry by also including "new geometries". *Geometric transformations*, in particular translations, were used to introduce the concepts of straight line and parallelism, as well as rotations were used to introduce the concept of angle and perpendicularity (as in Meray, 1874). In this period, new topics such as conic sections and analytical geometry made their first appearance in the school-curriculum.

Among those new texts we find a text by Julius Henrici and Peter Treutlein (1881-1883). It is not a book about practical geometry, even though it contains a part explicitly devoted to practical geometry. The extended use of *geometric transformations* is linked to movement and to geometric constructions: the authors state in the preface that it is important that the pupils concretely perform the transformations in the classroom using models.

In fact the importance of practical geometry is not so much as in its use for applications, but rather in its educational value; in particular in the interpretation of the term *practical* as "to be performed concretely", that is the meaning that it mainly acquired in the 1900s.

The text of Henrici and Treutlein is a text of synthetic geometry for the Gymnasium (for the Tertia, that is for pupils older than 13), with axioms and theorems. But it presents relevant differences from the Euclidean text. The authors state that they pursue "a logical order of the concepts rather than of the theorems". In particular the subject is developed according to the kind of transformations that are necessary for the proofs.

Here the use of geometric transformations is more relevant then in Méray. At first, figures with a centre are introduced, and the central symmetry is used to define parallel lines (p.23). To *draw* parallel lines the use of a set-square is also suggested (p. 21)

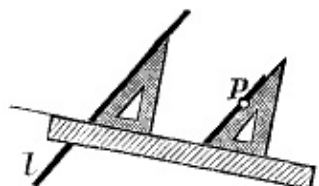

Fig. 39.

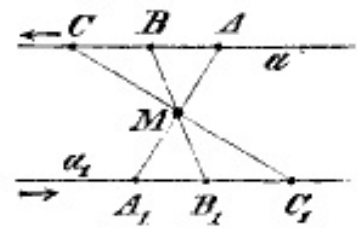

Fig. 43.

A *new kind of proof* can be performed using geometric transformations. Let's see how central symmetry is used to prove that the medians of a triangle are concurrent (p.36):

> If two of them, $AA_1$ and $BB_1$, meet in $S$, we can rotate the segment $CS$ around $B_1$ in $AB_2$ and around $A_1$ in $BA_2$, so that $CS$ is parallel and equal to $AB_2$ and $BA_2$. Also $AB_2$ and $BA_2$ can be superimposed and the centre of the rotation will be $S$. The three parallel lines $CS$, $AB_2$, $BA_2$ will cut equal segments on $AB$, on $AA_2$ and $BB_2$ [so $CS$ is a median].

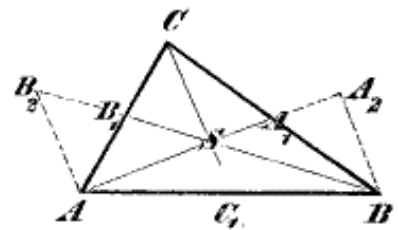

Fig. 77.

*A new kind of theorems* is also possible.

For instance, after the introduction of translations and rotations, the authors show how to find the transformations that bring two equal triangles to coincide (p.40).

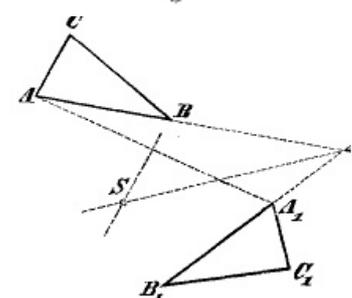

Fig. 84.

After the circle is introduced, a drawing with ruler and compass can be performed, and it is a relevant part. A *method* is given on how to face problems of constructions. At first, we find the usual problems that require the construction of perpendicular or parallel lines, as we can find in Clairaut, but we soon arrive at the more complicated problems, such as drawing tangents to a circle under certain conditions or drawing the common tangent to two given circles (p. 64).

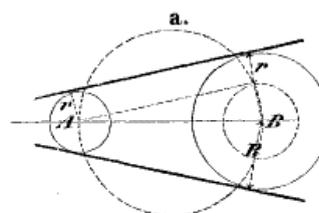
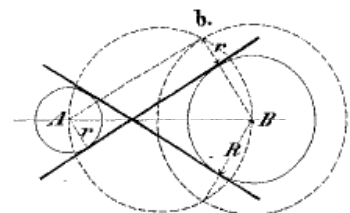

Fig. 121.

Areas are introduced with reference to Euclid's Book II (geometric algebra), to the proportions of the area of rectangles of equal bases, to equi-decomposition (seen also as equivalent transformation, p. 93).

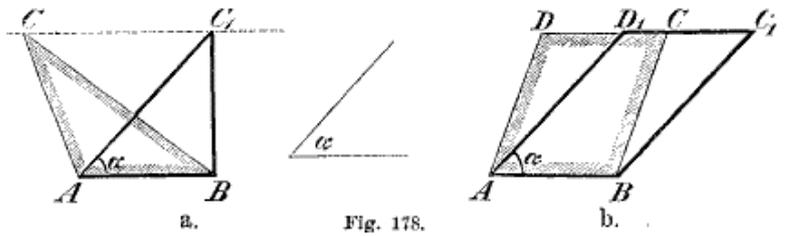

The theorems of Euclid and Pythagoras are thus proved with reference to equivalent transformations and not to proportions. This chapter is about the "transformations of figures", and contains exercises similar to those of Clairaut in his analogous chapter about "comparison of figures", but there are also problems of "divisions of figures" that are exactly like those of Fibonacci in his "division of fields", even if they have now lost they practical origin. For instance

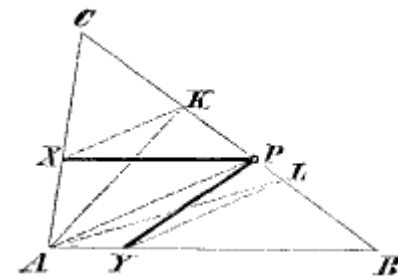

> To divide a triangle in *n* [3] parts from a point on a side (p. 96)

The solution is suggested by the figure.

In the appendix of the text by Henrici and Treutlein, the measurement of areas is introduced by giving the rule to calculate the area of a rectangle; this procedure is always that of Clairaut or Mocnik. We now find the *first numerical examples* of this text.

In the second volume similarity is introduced, and also elements of projective geometry. Besides elements of prospective representation, as we find in Mocnik, there is a more extended part containing topics such as harmonic points, pole and polar, projective correspondences, conic sections. Problems of measurement are taken up again and solved by algebraic means.

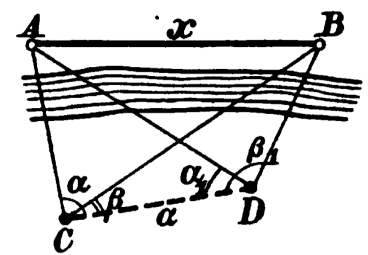

Problems of *practical geometry* are treated as an application of trigonometry. For instance, the problem that Mocnik solves with the help of similar representation and concrete measurement is solved by means of trigonometric formulas that lead to the value of *x* (p. 157):

$$x^2 = \left(\frac{a\sin\alpha_1}{\sin(\alpha + \alpha_1)}\right)^2 + \left(\frac{a\sin\beta_1}{\sin(\beta + \beta_1)}\right)^2 - 2 \cdot \frac{a\sin\alpha_1}{\sin(\alpha + \alpha_1)} \cdot \frac{a\sin\beta_1}{\sin(\beta + \beta_1)} \cdot \cos(\alpha - \beta)$$

Measuring is no longer the aim of teaching geometry, as it was for Fibonacci or Ramus; nor an "excuse", as it was for Clairaut (with the only exception of the applications of trigonometry). Numbers have a very small role. The aim of Henrici and Treulein is the comparison of figures by means of transformations. But the concrete references to the operations of transforming and the geometric constructions place this book in a new stream of synthetic "practical" geometry that evolved in the succeeding decades.

# 6. EXPERIMENTAL GEOMETRY

At the beginning of the 20[th] century an international reform movement started, which aimed at introducing modern topics, such as analytical geometry and calculus in secondary school, and also proposed new methodologies in mathematics teaching. As to geometry, a practical-intuitive approach was suggested. In many countries new syllabi were established. Borel (1905) wrote a text in accordance with the new French syllabi of 1905 (see Nabonnand, 2007). The various reprints of Henrici & Treutlein's text seemed to be suitable for the German syllabi of 1907 (Becker, 1994).

The *experimental geometry* of John Perry had a major influence at an international level. Perry (1901) emphasized the educational value of experimental procedures in the first approach to Euclidean geometry. He thought that a major part of elementary geometry had to be assumed as primitive and that the subject matter should be taught with reference to its utility as well as being interesting to the pupils (Howson, 1982; Barbin & Menghini, 2013).

The text of *practical geometry* by Joseph Harrison (1903) was written according to the proposals of Perry and even presented an appreciation of Perry himself in the preface.

Harrison states in the preface that many of the British schools and colleges are equipped with *laboratories* in which "experimental work involving quantitative measurements can be carried on as part of the ordinary school course", and it is "coming to be recognized that elementary mathematics should be taught in relation to such work".

The syllabus that the author cites at the end of the book (Science Subject I, Board of Education, South Kensington) states that

> This subject [practical plane and solid geometry] comprises the graphical representation of position and form and the graphical solution of problems. [...] in the elementary stage the main object of the instruction will be to familiarize the student with the fundamental properties of geometric figures and their applications [...] it is not intended that the student shall follow the Euclid's sequence [...]. In the Advanced stage [...] the subject will be developed in its applications [Engineering and physical group...].

We can see that the advanced stage doesn't propose Euclidean Geometry.

Harrison starts with a sort of definition of points, lines and surfaces which is a mixture of Euclid's definitions, Ramus' drawings and Mocnik's representation of points. The aim of drawing points is to explain how a point can be seen in a drawing.

Then he gives a detailed description of drawing instruments and their uses, and explains how to draw parallel and perpendicular lines using various kinds of instruments (ruler and set square, tee square and clinograph, p. 29).

Harrison gives a long explanation of how to measure angles, and also he defines the basic trigonometric functions and explains the

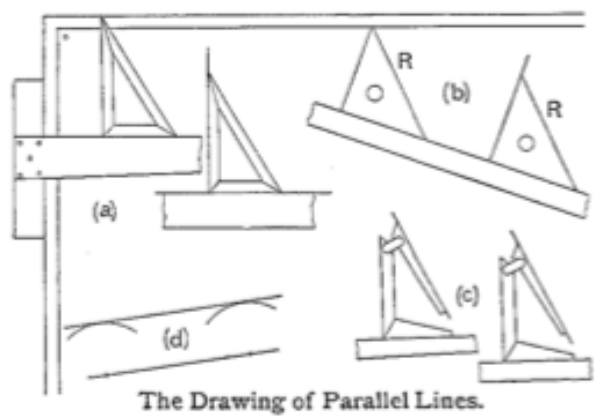

The Drawing of Parallel Lines.

trigonometric tables (so we can understand that the considered level is not a primary or beginning secondary school level).

One of the problems is, as was in Clairaut, that of determining the width of a river. The suggested method is different (please note that Harrison doesn't suggest to work in scale, so the solution is not very practical ...)

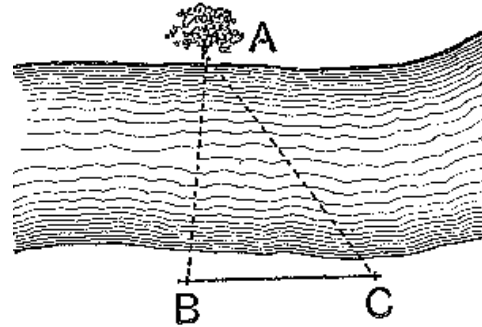

> The figure shows how the width of a river could be ascertained by a person on one bank. He might select and measure a base *BC*. Then note some conspicuous object *A* on the remote bank, and by a sextant or other instrument measure the two angles *ABC* and *ACB*. The triangle *ABD* could then be plotted, and the width measured (p. 61).

Defining a parallelogram Harrison states:

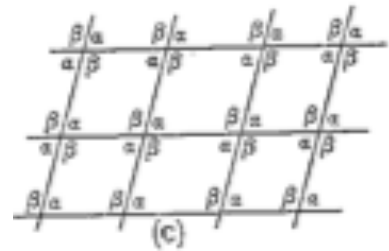

> A parallelogram is a quadrilateral in which opposite sides are parallel. Thus in fig. *c* the two systems of parallel sides give rise to a series of parallelograms. From the figure it is evident that opposite angles of a parallelogram are equal. And by measurement the student will easily satisfy himself that opposite sides are also equal (p. 42)

No further explanation is given.

Again through measurement or by paper cutting the student can verify the *Pythagoras theorem* (p. 47).

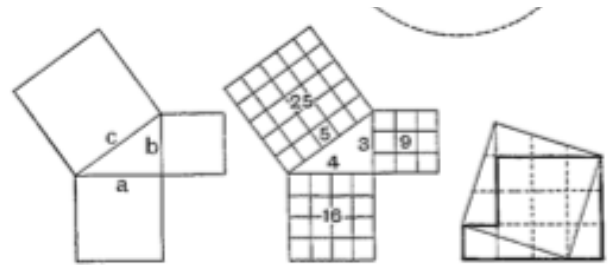

Harrison uses the same drawings of Fibonacci to show that the area of any plane figure may be measured by finding how many unit squares would be required to cover it (p. 77-79).

Another problem that we have already found in Clairaut is "To reduce a given rectangle *AC* to an equivalent rectangle on the base *BD*".

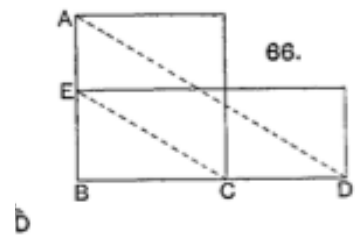

And its solution, different from Clairaut, is without explanations: "Join *DA*. Draw *CE* parallel to *DA*. Draw the rectangle ED, which is the one required" (p. 83).

The same kind of "proof" is given for other theorems, such as "The angle in a semicircle is a right angle":

> Verification: Draw any semicircle, diameter *AA*. Set out several angles *ABA* in the semicircle. Verify that in all cases the angle *ABA* is 90°.

Harrison introduces also vectors, their representation and sum and difference. Practical geometry in space is about methods of representing solids, mainly by projections on the coordinate planes. These exercises are definitely not elementary technical drawing exercises. They require, for instance, the use of trigonometry to calculate angles and measures.

Not all English books on practical or experimental geometry are of the same kind; in the text by Godfrey and Siddons (1903) the part devoted to measurement is smaller than Harrison's, moreover and it presents exercises which require an idea or a form of explanation (Fujita et al., 2004). Anyway, the main focus of the practical geometry in that period is on measuring or calculating measures, with some technical drawing; it is quite obvious the we would find objections to this methodology:

> Most recent text-books of geometry contain an introduction on practical geometry. While presupposing a short preliminary course of this nature, we have preferred to leave it to the teacher to devise himself. In this direction we think that the recent reforms have gone too far, and we feel sure that, as regards secondary schools, it will be necessary to retrace our steps. Too much time spent on experimental and graphical work is wearisome and of little value to intelligent pupils. They can appreciate the logical training of theoretical geometry, while experiments of far greater interest can be made in the physical and chemical laboratories (preface, Davison, 1907).

## 7. INTUITIVE-INTRODUCTORY GEOMETRY

Notwithstanding these criticisms, experimental geometry has influenced, from the beginning of the 1900s, texts of *intuitive or introductory geometry* devoted to the lower school grades. These texts were influenced not only by experimental geometry, but also by the whole stream of texts of practical geometry (Menghini, 2010). An earlier intuitive - experimental approach was considered a good aid for students to overcome the difficulties caused by the logical deduction of Euclid's textbook.

So we find *free-hand drawing* (Veronese, 1901) and *drawing instruments* (Borel, 1905), and parallel lines introduced by means of central symmetry or translation. *Geometric transformations* were considered not only useful to introduce simple concepts, but also to compare segments and figures, as they were carried out practically. The equality of triangles was established through the possibility of constructing them (as in Clairaut), and *measurement* and the use of *numbers* became widespread.

We also find the use of *paper folding* and *concrete materials.* For instance: two successive foldings of a piece of paper are used to introduce the concept of perpendicularity (Frattini, 1901; Borel, 1905); in order to prove that *the diagonals of a parallelogram bisect each other,* a parallelogram is cut out from a piece of paper, filling again the resulting empty space after a rotation of 180 degrees (Frattini, 1901); to know the sum of the angles of a triangle, the corners of a triangle drawn on paper are cut and placed next to each other to check that they form a straight angle (Amaldi, 1941).

Emma Castelnuovo, who was also largely inspired by the approach via problems of Clairaut, uses simple tools, such as a folding meter, to show how to modify a quadrilateral into a different one and to analyze the limit situations (Castelnuovo, 1946).

# 8. DEVELOPMENTS OF PRACTICAL GEOMETRY IN THE 20^{TH} CENTURY

## 8.1. Axioms for measure

The textbook *Basic Geometry*, by George David Birkhoff and Ralph Beatley, was published in an experimental edition in 1933, and officially edited in 1941. The text starts with many exercises that require the practical use of a scale and a protractor: measuring of lengths and angles, and verifying rules about their addition and subtraction. Numbers are very important from the beginning. It would seem a text about practical geometry, but *measuring is not an aim*. The aim is to help pupils to accept the *Principles*. The first is about:

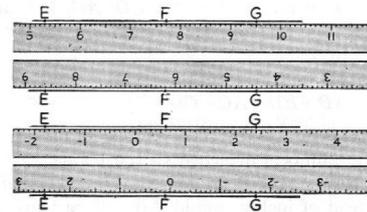

Fig. 2 shows what happened when four different pupils attempted to measure the distance *EG* on a certain straight line. They all used the same scale to measure this distance. Did they all get the same result?

Fig. 2

> *Line measure*: the points of any straight line can be numbered so that number differences measure distances (p. 40).

And the third about:

> *Angle measure*: All half lines having the same end point can be numbered so that number differences measure angles (p. 47).

Two more axioms state that "there is one and only one straight line through two given points" and that "all straight angles have the same measure", while the fifth axiom corresponds to a criterion of *similarity of triangles*. An example of the application of the axioms is the proof of the theorem about the *sum of the angles of a triangle*:

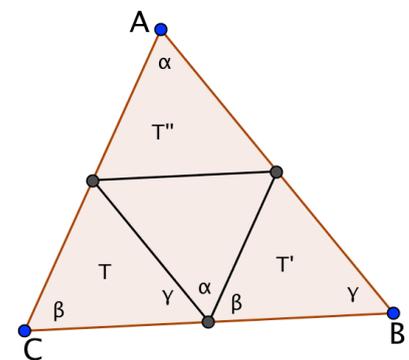

Principle 1 guarantees that the sides of the triangle *ABC* can be divided by two, while Principle 5 guarantees the similarity of *T*, *T'*, *T''* to *ABC*, and thus the similarity of *ABC* and the internal triangle. So ... (see figure).

The ideas of Birkhoff slowly found their way into American teaching, and inspired many other textbooks on geometry, such as those produced by the School Mathematics Study Group.

## 8.2. Old applications and new theorems

In the period of the New Math reforms in the 60s, the British *School Mathematics Project* (1970) still presented elements of practical geometry. A chapter was devoted to the measuring of an angle, to the use of the protractor, and to applications such as:

> To measure the height of a hill, land-surveyors measure the angles *a* and *b* from two different points *A* and *B* which are at the same height, and whose height and distances are known. If *AB* = 1000 m, *a* = 19° and *b* = 36°, draw a diagram in a scale 1 cm to 200 m and find *DC*. If *A* and *B* are at 200 m height, what is the height of the hill? (p. 82)

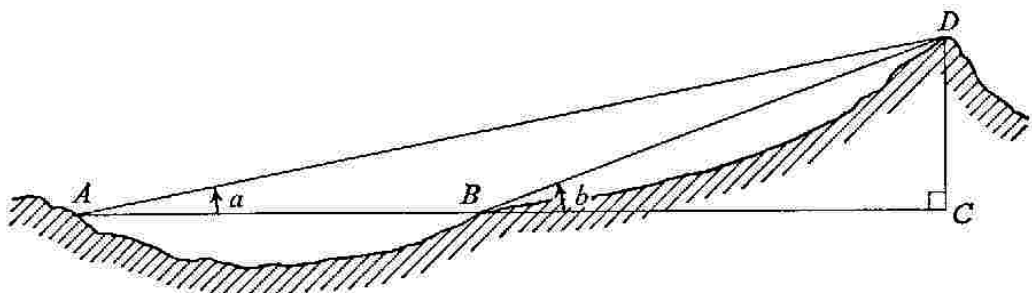

Pupils are supposed to find the answer by measuring, as was seen in exercises of Clairaut or Mocnik.

After chapters about polygons and polyhedra, approximations of areas by tessellation, observation of the symmetries of a figure, and paper folding, a chapter is again devoted explicitly to land-surveying. The aim is not only to apply what has been studied in the previous chapters, but also to introduce new concepts. Here the old applications of practical geometry are the basis for exercises in the new practical geometry. There is also the suggestion to do real exercises in the open air.

These exercises mainly require to draw on paper a piece of land by measuring (or given) certain distances and angles using similarity (p. 298); *irregular* forms (p. 300) can be drawn by using the way proposed by Mocnik.

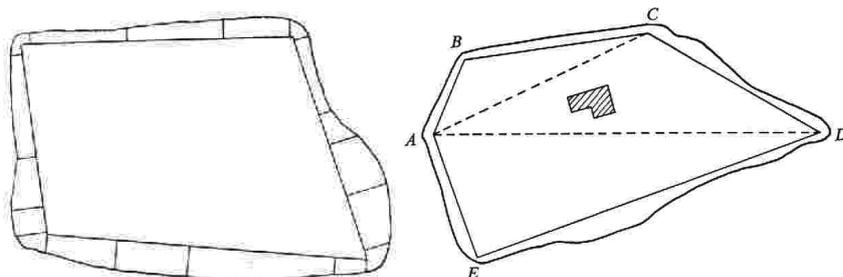

So we find many of the topics that we have met from the beginning of practical geometry.

The book also presents *geometric transformations* and *new theorems* based on their *composition* (as we have seen in Henrici & Treutlein) and on their *invariants*.

**CONCLUSION**

For Fibonacci instruments and their use to measure distances were taken for granted. With the exception of a short part about the unit measures used at that time, his text is mainly concerned with rules for calculating areas and with theorems that can help in comparing the areas. The rules arise mostly from a generalization of numerical examples, and numerical values are given to the considered segments in order to verify theorems. These simplifications are an aid for the reader, who at the same time is challenged with concrete problems that require finding a distance (using similarity), or with more abstract problems that require the division of a field given certain conditions. For the solution to these problems only an initial hint to allowed constructions is given.

With Petrus Ramus we find on the one hand an increase in the *activity* of measuring, through the explicit presentation of measuring and drawing *instruments*, on the other hand this activity is postponed in order to *present* geometric figures, their properties, and to "joke" with them.

We can see that the "puzzles" of the practical geometry of the Middle Ages and the *jokes* of Ramus are characteristic of a *freedom* that is not in the Euclidean tradition. Euclid is absolutely not "playful".

Clairaut focuses again on activity, which is now aimed at *constructing* the geometric elements to be measured. As in Fibonacci, not only is measuring per se is not important, but also the starting problems are not as important as the process of solving them. The role of the process of constructing is principally an *educational* one.

The educational role of constructions, of concrete activities, and of observation can also be found in the text of Mocnik, whose aim is again practical. Finding distances and areas is still a goal, but the use of measuring instruments is in the background, while theorems acquire importance.

The text of Henrici and Treutlein cannot be considered exactly in the same line of thinking, but it takes from practical geometry examples, problems and an idea of concrete activity, which is then applied to new objects as the geometric transformations. Here we can see how practical geometry suggested a methodology that also influenced the more theoretical study at the Gymnasium.

In the *practical geometry* of the 1900s we find something similar to Fibonacci. If Fibonacci verified the theorems by substituting lengths with numbers, Harrison verifies them by measuring. The use of instruments and the activity of measuring are central to Harrison; but we find nothing that can be considered an intellectual challenge to the pupil.

The *intuitive geometry* for the lower grades can be considered as a middle ground: those activities of measuring, of concretely performing transformations, of using instruments, can be seen as an introductory activity *before* starting with Euclidean deduction.

Another way of combining practical and deductive geometry is proposed in the book of Birkhoff and Beatley, where measuring must help pupils to understand the role and the meaning of axioms for line and angle measure. This book shows that a *deductive* geometry based on measure needs new axioms.

Avoiding axioms, the SMP proposes a curriculum that encompasses many of the topics that we have seen in the development of practical geometry: from the concrete problems of land-surveying, to drawing, to the *new theorems and proofs* based on geometric transformations.

The proposals of SMP in the 60s show that practical geometry can have a role that goes beyond an introductory or intuitive presentation. As we understand from the presented development, a practical-intuitive teaching of geometry requires an active role of the pupil, where reasoning is supported by action and by construction. We have also seen examples of rational teaching of geometry supported by practical elements. Even if the challenging problems proposed by the practical geometry of the late Middle Ages entered only occasionally in the successive geometry textbooks (with the exception of the *division of fields*), we can recognize the heritage of this geometry in a certain freedom in choosing problems and methods (including *intuition*). However the main methodological heritage lies in the shift of the meaning of "practical" from "useful for applications" to "to be performed concretely". This shift happened towards the end of the Middle Ages even in texts whose aim was very concrete, through the *description* of the activity of measuring or of constructing, and the use of instruments.

Practical geometry can thus suggest alternative methodologies, not only for an intuitive introduction *before* rational geometry, not only for examples of useful applications, but also for a different way of teaching based on a conception of knowledge which includes the *practical process of knowing*.